\documentclass[10pt]{amsart}
\usepackage[cp1251]{inputenc}
\usepackage[english,english]{babel}
\usepackage{amsmath}
\usepackage{amssymb}
\usepackage{amsfonts}

\newcounter{subsect1}

\newfont{\gotikai}{eufm10}

\def\v{\varepsilon}

\def\d{\delta}
\def\a{\alpha}
\def\b{\beta}
\def\g{\gamma}
\def\l{\lambda}
\def\m{\mu}
\def\n{\nu}
\def\t{\tau}
\def\x{\xi}
\def\z{\zeta}

\def\h{\theta}

\newcommand{\skk}[1]{\left\{ #1 \right\}}

\def\DD{P\ r\ o\ o\ f}

\def\R{\mathbb R}
\def\Z{\mathbb Z}

\def\x{\xi}

\begin{document}

\title[Local theorems for arithmetic renewal
processes]{LOCAL THEOREMS FOR ARITHMETIC COMPOUND RENEWAL PROCESSES, WHEN CRAMER'S CONDITION HOLDS}
\author{{A.A. MOGULSKII}}%
\address{Anatolii Alfredovich Mogulskii  
\newline\hphantom{iii} Sobolev Institute of Mathematics,
\newline\hphantom{iii} pr. Koptyuga, 4,
\newline\hphantom{iii} 630090, Novosibirsk, Russia%
\newline\hphantom{iii} Novosibirsk State University,
\newline\hphantom{iii} 1 Pirogova Str.,
\newline\hphantom{iii} 630090, Novosibirsk, Russia}%
\email{mogul@math.nsc.ru}%
\thanks{\sc  Mogulskii A.A. Local theorems for arithmetic compound renewal processes, when Cramer's condition holds.}
\thanks{\copyright \ 2018 Mogulskii A.A.}
\thanks{\rm The research was supported by Russian Science Foundation (grant No. 18-11-00129).}

\maketitle {\small
\begin{quote}
\noindent{\sc Abstract. }
We continue the study of the compound renewal processes (c.r.p.), where the moment Cramer's
condition holds (see \cite{I1}---\cite{I.10}, where the study of c.r.p. was started).
In the paper arithmetic c.r.p. $Z(n)$ are studied. In such processes random vector $\xi = (\tau,\zeta)$
has the arithmetic distribution, where $\tau >0 $ defines the distance between jumps, $\zeta$
defines the values of jumps. For this processes the fine asymptotics in the local limit theorem
for probabilities  ${\bf P}(Z(n)=x)$ has been obtained in Cramer's deviation region of $x\in \Z$.
In \cite{I6}---\cite{I.10}  the similar problem has been solved for non-lattice c.r.p., when
the vector $\x=(\t,\z)$ has the non-lattice distribution. \medskip

 \end{quote}
}

\setcounter{section}{1}
\setcounter{equation}{0}

                       \begin{center}
                       {\bf \S~1. Introduction. Formulation of the problem.                        }
                    \end{center}


     Consider the sequence
 \begin{equation}{\label{1.1}}
  \{\x_k=(\t_k,\z_k)\}_{k=0}^\infty
\end{equation}
of random vectors, such that $\t_0=0$, $\z_0=0$, $\t_1\ge 0$, $\t_k>0$ when $k\ge 2$.
 Denote
 $$
T_n:=\sum_{k=0}^n\t_k,~~~Z_n:=\sum_{k=0}^n\z_k,
 $$
 such that $T_0=Z_0=0$, and where $n\ge 0$.
Let $\n(0):=0$,  and for $t > 0$
\begin{equation}{\label{1.2}}
\n(t)=\max\{k\ge 0:~T_{k}<t\}.
\end{equation}

{\it Compound renewal process} (c.r.p.)
$Z(t);~t\ge 0$, is defined by the equality (see, for example, \cite{I1},
\cite{I6})
\begin{equation}{\label{1.3}}
  Z(t):=Z_{\n(t)},
\end{equation}
when  (\ref{1.1}) is the sequence of independent random vectors,
where vectors $\x_k=(\t_k,\z_k)$ when $k\ge 2$ have the same distribution as the vector
$\x=(\t,\z)$.
Standard conventional model of c.r.p. assumes that the joint distribution of the moment $\t_1$
of the first jump and its size $\z_1$ differs em general from the
joint distribution of $(\t,\z)$ (see, for example, \cite{I1},
\cite{I8},
\cite{I9}).
It is implemented, for example, for c.r.p. with stationary increments.

If $(\t_1,\z_1)=_d(\t,\z)$,
 then the process $Z(t)$ we call {\it homogeneous c.r.p.};
 and {\it nonhomogeneous } otherwise.
 The process  $\n(t)$  (it corresponds to the case when $\z_1=\z_2=\cdots =1$)
we will call {\it simple} renewal process.

Thus, the distributions of two random vectors
$\x_1=(\t_1,\z_1)$, $\x=(\t,\z)$
{\it  determine completely } the distribution of c.r.p. $Z(t);~t\ge 0$.
Under the event
$$
     \{T_k<t\le T_{k+1}\}~~~~\mbox{when}~~~k\ge 0
$$
the following equalities hold true: $\n(t)=k$,  $Z(t)=Z_k$. Therefore,
the step processes
 $\n(t)$, $Z(t)$
are left continuous when $t>0$.

From the strong law of large numbers for the sums $(T_n,~Z_n)$
it is follows that, when $t\to \infty$
$$
  \frac{\n(t)}{t}\to \frac{1}{{\bf E}\tau},~~~
  \frac{Z(t)}{t}\to a:=\frac{{\bf E}\z}{{\bf E}\tau}
$$
with probability 1.

In this work, as in \cite{I6}, \cite{I.10}, we will suppose that the random vectors
$\x_1=(\t_1,\z_1)$, $\x=(\t,\z)$, which define c.r.p. $Z(t)$, satisfy the following
moment Cramer's condition

\begin{enumerate}
 \item[]
$[{\bf C}_0]$. (Cramer's condition) {\it ${\bf E}e^{\d|\x_1|}<\infty$,
${\bf E}e^{\d|\x|}<\infty$
for some $\d>0$}.
\end{enumerate}

In papers \cite{I6}, \cite{I7} integro-local theorems for the
process $Z(t)$, when $t\to \infty$, were obtained in the case, when the random vector
$\x=(\t,\z)$, ``determining" the process $Z(t)$, is not lattice distributed.
In this paper the analogous problem is being solved for
{\it arithmetic} c.r.p. $Z(t)$, it means that the values of random vectors $\x_1$, $\x$
belong to the integer lattice $\Z^2$.  More precisely, we suppose that
${\bf P}(\x_1\in \Z^2)=1$, and the following more strong {\it arithmetic} condition
for the vector $\x$ in terms of its characteristic function holds true
$$
  f(u):={\bf E}e^{i\langle u,\x\rangle},~~~u\in {\mathbb R}^2,
$$
where $\langle \a,\b\rangle$ is scalar product in $\R^2$.

\begin{enumerate}
 \item[]
$[{\bf Z}]$. (Arithmetic condition) {\it
For any $u\in \Z^2$ the equality $f(2\pi u)=1$ holds true, and
the inequality $|f(2\pi u)|<1$ holds for any
$u \in \R^2\setminus \Z^2$.}
\end{enumerate}

Further, {\it  we always assume that the moment conditions $[{\bf C}_0]$ and the structural arithmetic
condition $[{\bf Z}]$ hold; thus, in order to avoid repetitions the fulfillment of conditions $[{\bf C}_0]$, $[{\bf Z}]$ will not be recalled}.

Note that, if the condition $[{\bf Z}]$ holds for the random vector $\x=(\t,\z)$,
then this vector cannot be degenerate in space $\R^2$, that is, with probability
1 it belongs to a line
$$
  L_{n,c}:=\{(u,v)\in \R^2:~~n_1u+n_2v=c\},
$$
which is defined by a normal vector $n=(n_1,n_2)$ and constant $c$. Indeed,
if ${\bf P}(\x\in L_{n,c})=1$,
then there exists $r$ such that $u=(u_1,u_2):=(rn_1,rn_2)\not \in \Z^2$ and the equality
$$
  |f(2\pi u)|=|{\bf E}e^{i2\pi rc}|=1,
$$
holds true, which contradicts with the condition $[{\bf Z}]$.

Since the arithmetic c.r.p. $Z(t)$ changes its values only at integer moments
 $t=1,2,\cdots,$ and it is left continuous, then for any integer
 $t> 0$ the equality $Z(t)=Z([t]+1)$ holds true, where, as usual,
 $[t]$ denotes the integer part of nonnegative number $t$.
 Therefore we will consider only integer values of $t=n\in \{ 0,1,2,\cdots\}$,
 that is we will study the random sequences
 $$
   Z(0),~Z(1),~Z(2),\cdots.
 $$
 It simplifies notations and it does not lead to a loss of generality.
 In the present work we established the local theorem for arithmetic c.r.p. $Z(n)$,
 where the explicit asymptotics for the probability
 $$
   {\bf P}(Z(n)=x)\sim ?~~~\mbox{when}~~~n\to \infty,
 $$
 was found for the sequence $x=x_n\in \Z$, such that the point $\a:=\frac{x}{n}$
 belongs to some fixed compact set $K\subset \R$  (see Theorem~2.1 below).

 The proof of the local theorem for arithmetic c.r.p. $Z(n)$ is based on
 local limit theorem for the renewal function (measure)
$$
H(B):=\sum_{n=0}^\infty {\bf P}(S_n\in B),~~~B\subset \R^2,
$$
where
the explicit asymptotics was studied for the set $B=B_{t,x}:=\{t\}\times\{x\}$
$$
  H(B_{t,x})\sim ?,
  $$
for the sequence $(t,x)=(t_n,x_n)\in \Z^2$ such that the point
$(\theta,\a):=\frac{1}{n}(t,x)$ belongs to some fixed compact set
$K\subset \R^2$ (see Theporem~2.2).


The used methods follow the papers \cite{I6},
\cite{I7}, where under the moment Cramer's condition $[{\bf C}_0]$
the {\it integro-local limit theorem} was obtained for ``non-lattice"
c.r.p. $Z(t)$ and for corresponding renewal measure $H(B)$, that is, the exact asymptotics were found
$$
  {\bf P}(Z(T)\in \Delta[x))\sim ?,~~~H(\d[t)\times\Delta[x))\sim ?,
$$
where $\Delta[x):=[x,~x+\Delta)$, $\d[t):=[t,~t+\d)$ are
semi-intervals of length $\Delta$, $\d$, respectively.
Therefore Theorems 2.1 and 2.2 in this work complement the results of
the paper
\cite{I6} (Theorem 1.1 and 3.1, respectively).

Historical reviews of related results can be found in
\cite{I1}, \cite{I6}. We add only the references to the works
\cite{I.8}---\cite{I.10}, where integro-local limit theorems
for {\it multidimensional} non-lattice c.r.p. were obtained under the
Cramer's condition $[{\bf C}_0]$.

The rest of the work contains on \S\S~2,3. In \S~2 the main results
are formulated. A number of functions will be involved in these formulations.
In order to understand the essence of established results
it is desirable to know the meaning and properties of these functions.
Therefore, following the paper \cite{I6},
in \S~2 we start by introducing the necessary concepts and notations
which are supplied with the necessary explanations and references.
The proofs of the main statements are in \S~3.

\setcounter{section}{2}
\setcounter{equation}{0}

                       \begin{center}
                       {\bf \S~2. Main results: local theorems for
                       renewal functions and compound renewal processes.
                        }
                       \end{center}

{\bf 2.1. First and second deviations functions.}
Basically, the notations follow the paper \cite{I6}.
Let
$$
\psi(\l,\m):={\bf E}e^{\l\t+\m\z},~~~
\psi_1(\l,\m):={\bf E}e^{\l\t_1+\m\z_1},
$$
$$
  A(\l,\m):=\ln \psi(\l,\m),
  ~~~A_1(\l,\m):=\ln \psi_1(\l,\m),~~~(\l,\m)\in \R^2;
$$
$$
  \mathcal{A}:=\{(\l,\m):~\psi(\l,\m)<\infty\},~~~
  \mathcal{A}_1:=\{(\l,\m):~\psi_1(\l,\m)<\infty\}.
$$
It is clear that according the condition $[{\bf C}_0]$
the {\it interiors} $(\mathcal{A})$, $(\mathcal{A}_1)$
of the sets $\mathcal{A}$, $\mathcal{A}_1$
contain the point $(\l,\m)=(0,0)$, and they are the areas of analyticity of the
functions $A(\l,\m)$, $A_1(\l,\m)$, respectively.
In the proof of the main results of this work (as in \cite{I6}, \cite{I7})
the known integro-local theorems for the sums
$S_n=(T_n,Z_n)$ (see, for example, \cite{I1}, \S~2.9) are used.
Where the important role is performed by the so-called
{\it (first) deviation function }
\begin{equation}{\label{2.1}}
\Lambda(\h,\a):=\sup_{(\l,\m)}\{\l \h+\m\a-A(\l,\m)\},
 \end{equation}
which correspond to the random vector $\x=(\t,\z)$. It is the Legendre transform
of the convex function $A(\l,\m)$;
because why the function $\Lambda(\h,\a)$ is also convex and continuous from below.

We agree on the following notations. Further, for a given function
$F=F(u,v)$ of two variables $u$ and $v$, the subscripts
$(1)$ and $(2)$ will denote the derivatives on the first and second argument
respectively, for example:
$$
F'_{(1)}(u,v)=\frac{\partial}{\partial u} F(u,v)|,~~~
F''_{(2,1)}(u,v)=
\frac{\partial}{\partial u}
\frac{\partial}{\partial v}F(u,v)|.
$$
$F'=F'(u,v)$ and $F''=F''(u,v)$ will denote the vector
$$
F'=F'(u,v)=(F_{(1)}'(u,v),F_{(2)}'(u,v))
$$
and the matrix
$$
F''=F''(u,v)=\|F_{(i,j)}''(u,v)\|_{i,j=1,2}.
$$
$|F''|$ denotes the determinant of the matrix $F''$. For the matrix
$M=||M_{i,j}||_{i,j=1,2}$ through
$$
  (\h,\a)M(\h,\a)^\top:=\h^2M_{1,1}+\h\a(M_{1,2}+M_{2,1})+\a^2M_{2,2}
$$
we denote the corresponding quadratic form (where $^\top$ means the transpose,
it transforms a row into a column).

Together with the sets $\mathcal{A}_1$ and $\mathcal{A}$ we will need
the area $\mathcal{L}$ of analyticity of the function $\Lambda(\h,\a)$.
This area contain the points $(\h,\a)$ for which
the system of equations
 for the points coordinates $(\l,\m)$ (where supremum in (\ref{2.1}) is attained)
\begin{equation}{\label{2.2}}
\left\{
  \begin{array}{ll}
   A'_{(1)}(\l,\m)=\h,\\
    A'_{(2)}(\l,\m)=\a,
  \end{array}
\right.
 \end{equation}
 has the root $(\l(\h,\a), \m(\h,\a))$ belonging $(\mathcal{A})$ such that
$$
   \mathcal{L}=\{A'(\l,\m):~(\l,\m)\in (\mathcal{A})\}.
$$
Since the function $A(\l,\m)$ is strictly convex on $(\mathcal{A})$, then
this solution is always unique. It means that the conditions
$$
  (\h,\a)\in \mathcal{L}~~~\mbox{and}~~~(\l(\h,\a), \m(\h,\a))\in (\mathcal{A})
$$
are equivalent. Herein the mappings
$$
  (\h,\a)=A'(\l,\m),~~~(\l,\m)=\Lambda'(\h,\a)
$$
are mutually inverse; they are bijections between the sets
$(\mathcal{L})$ and $\mathcal{L}$.
It is clear that the point $(\h,\a)=(a_\t,a_\z)$, where $a_\t:={\bf E}\t$, $a_\z:={\bf E}\z$,
always belongs to $\mathcal{L}$; and
$(\l(a_\t,a_\z), \m(a_\t,a_\z)=(0,0)\in (\mathcal{A})$.
It is known (see, for example, \cite{I6}),  that for $(\h,\a)\in \mathcal{L}$
the following is true
$$
 \Lambda'(\h,\a)=(\Lambda'_{(1)}(\h,\a),~\Lambda'_{(2)}(\h,\a))=
 (\l(\h,\a),~\m(\h,\a)),
$$
where the pair of the functions $\l(\h,\a),~\m(\h,\a)$ is the unique solution of
the system (\ref{2.2}).

If $\t$ and $\z$ are independent, then
$$
  A(\l,\m)=A_\t(\l)+A_\z(\m),
$$
where
$$
 A_\t(\l):=\ln {\bf E}e^{\l\t},~~~A_\z(\m):=\ln {\bf E}e^{\m\z}.
$$
Thus $A'_{(1)}(\l,\m)=A'_\t(\l)$ does not depend on $\m$,
$A'_{(2)}(\l,\m)=A'_\z(\m)$ does not depend on $\l$, and
the areas $(\mathcal{A})$ and $\mathcal{L}$ are rectangulars.

With the deviations function $\Lambda(\h,\a)$ we need also a so-called
{\it second deviations function $D(\h,\a)$}, which is defined by
 \begin{equation}{\label{2.3}}
D(\h,\a):=\sup_{(\l,\m)\in \mathcal{A}^{\le 0}}\{\l\h+\m\a\}=
\sup_{(\l,\m)\in \partial\mathcal{A}^{\le 0}}\{\l\h+\m\a\},
 \end{equation}
where $\mathcal{A}^{\le 0}:=\{(\l,\m):~A(\l,\m)\le 0\}$, and $\partial B$ is
the boundary of $B$ (see \S~2.9 from \cite{I1}).
The properties of the function $D(\h,\a)$ are fully studied in
\cite{I1}, \cite{I6}.  It is {\it convex, semi-additive, linear on any ray from
the center $0$, continuous from below}. Moreover, for all points
$(\h,\a)$ except the boundary points of the set $\mathcal{D}^{<\infty}$
of finiteness of the function $D(\h,\a)$ the following equality holds true
$$
D(\h,\a)=\inf_{r>0}r\Lambda\left(\frac{\h}{r},\frac{\a}{r}\right).
$$

Representation (\ref{2.3}) allows us to find one more useful characterization of
the function $D(\h,\a)$. We note first that due to linearity of the function
$D$ the following holds for $\h>0$
 \begin{equation}{\label{2.4}}
D(\h,\a)=
\h D\left(1, \frac{\a}{\h}\right),
 \end{equation}
 that is the function of two variable $D(\h,\a)$
 is completely determined by the function of one variable
 $$
  D(\a):=D(1,\a).
 $$
 Furthermore the main results of the paper are formulated in terms of
 this function.
Due to (\ref{2.3}) we have
 \begin{equation}{\label{2.5}}
D(\a)=
\sup_{(\l,\m)\in \partial\mathcal{A}^{\le 0}}\{\m\a+\l\}.
 \end{equation}
In order to describe the boundary $\partial\mathcal{A}^{\le 0}$ we consider the section
$\mathcal{A}_\m$ of the set
$\mathcal{A}$ on the level $\m$:
$$
 \mathcal{A}_\m:=\{\l:~(\l,\m)\in \mathcal{A}\}
$$
and let
$$
  \m^+:=\sup\{\m:~\mathcal{A}_\m\not = \emptyset\},~~~
\m^-:=\inf\{\m:~\mathcal{A}_\m\not = \emptyset\}.
$$
For any $\m\in (\m^-,\m^+)$
the function $A(\l,\m)$ is increasing on $\l$ within interval $(-\infty, \infty)$ and,
the following values are defined
$$
A(\m):=-\sup \{\l:~A(\l,\m)\le 0\}, ~~~A^\infty(\m):=-\sup \{\l:~A(\l,\m)<\infty\}.
$$
It is obvious that $(-A(\m),\m)\in \partial\mathcal{A}^{\le 0}$,
$(-A^\infty(\m),\m)\in \partial\mathcal{A}$,
it means that $A(\m)$, $A^\infty(\m)$ are finite and convex in the interval
$(\m^-,\m^+)$. Define the function $A(\m)$ out of the interval $(\m^-,\m^+)$ preserving
convexity and continuity from below by the following
$A(\m):=\infty$ when $\m\not \in [\m^-,\m^+]$  and
$$
 A(\m^+):=\lim_{\m\uparrow\m^+}A(\m),~~~A(\m^-):=
 \lim_{\m\downarrow\m^-}A(\m).
$$
Thus, the formula
(\ref{2.5}) can be written as
 \begin{equation}{\label{2.6}}
D(\a)=
\sup_{\m}\{\m\a-A(\m)\},
 \end{equation}
 it means that $D(\a)$ is {\it Legendre transform} on $A(\m)$.
 In \cite{I5} it was shown that the function
 $A(\m)$ possesses many properties of the logarithm of the Laplace transform over some distribution.
 It was noted that it is always convex and continuous from below on $\mathbb{R}$.
 Moreover, $A(\m)\to \infty$ when $|\m|\to \infty$, and if random variable $\z$ is not trivial,
 $$
 A'(0)=a:=\frac{a_\z}{a_\t}=\lim_{t\to \infty}{\bf E}\frac{Z(t)}{t},
 $$
$$
   A''(0)= \sigma^2:=\frac{1}{a_\t}{\bf E}(\z-a\t)^2=
  \lim_{t\to \infty}{\bf D}\frac{Z(t)}{t}.
$$
All the above function properties $A(\m)$ can be obtained from
(\ref{2.6}).
Obviously that the inequality $A^{\infty}(\m)\le A(\m)$ always holds true.
Since $(0,0)\in (\mathcal{A})$,
$A(0,0)=0$, then in the neighborhood of the point $\m=0$ the inequality
$A^{\infty}(\m) < A(\m)$ holds true.
Set
$$
  \m_-:=\max\{\m<0:~A(\m)=A^{\infty}(\m)\},~~~
  \m_+:=\min\{\m>0:~A(\m)=A^{\infty}(\m)\}.
$$
In the interval
 $\m\in (\m_-,\m_+)$
the following holds true $$
A^{\infty}(\m)<A(\m),~~~
(-A(\m), \m)\in (\mathcal{A}),$$
and  $\l=-A(\m)$ is the unique solution of the equation
\begin{equation}{\label{2.7}}
  A(\l,\m)=0.
\end{equation}
Thus, according the implicit function theorem $-A(\m)$ is the analytic function.
Let $\m(\a)$ be the point where the $\sup$ is attained
in (\ref{2.5}). If $A(\m)$ is differentiable in the point $\m(\a)$, then the value
$\m(\a)$ is the solution of equation
$$
  A'(\m)=\a.
$$
Since the function $A'(\m)$ increases monotonically, then there exists a unique solution
of the equation
$$
  \m(\a)=(A')^{(-1)}(\a),
$$
which is inverse function of $A'(\m)$. Let
$$
\a_{\pm}:=A'(\m_{\pm}\mp 0).
$$
  Then it is clear that
  the function $\m(\a)$ will be analytic in the interval $(\a_-,\a_+)$. Further, in
  \cite{I5} the following relation was established
  \begin{equation}{\label{2.8}}
D(\a) = \a\m(\a)-A(\m(\a))=\int_a^\a\m(v)dv.
 \end{equation}
From (\ref{2.8}) it follows that the function $D(\a)$ is analytic function on $(\a_-,\a_+)$.
Further, if we let
 \begin{equation}{\label{2.9}}
\l(\a):=-A(\m(\a)),
  \end{equation}
then from (\ref{2.7}) and (\ref{2.8}) we obtain
 $$
A(\l(\a),\m(\a))=0,~~~D(\a)=\l(\a)+ \a\m(\a)~~~\mbox{when}~~~\a\in (\a_-,\a_+).
$$

  From the above and from (\ref{2.4}) we obtain that the function $D(\h,\a)$
  is finite in the cone
  $$
    \mathcal{D}^{<\infty}:=\skk{(\h,\a):~\h>0,~\m\left(\frac{\a}{\h}\right)\in (\m^-,\m^+)},
  $$
  and it is analytic in the cone
  $$
 \mathcal{D}:=\skk{(\h,\a):~\h>0,~\frac{\a}{\h}\in (\a_-,\a_+)}.
$$
 In \cite{I6}, \cite{I7} it was shown that when
 $(\h,\a)\in \mathcal{D}$ then the following holds true
 $$
D'(\h,\a)=(D_{(1)}'(\h,\a),D_{(2)}'(\h,\a))=
\left(\l\left(\frac{\a}{\h}\right), \m\left(\frac{\a}{\h}\right)\right).
$$

{\bf 2.2. Local theorem for arithmetic compound renewal process.}
We will study the asymptotics of the probability
   \begin{equation}{\label{2.10}}
{\bf P}(Z(n)=x)~~~~\mbox{when}~~~x\in \Z,
 \end{equation}
 in the region of normalized deviations
   \begin{equation}{\label{2.11}}
\a:=\frac{x}{n}\in (\a_-,\a_+), ~~~x\in \Z.
 \end{equation}
 But in \cite{I6} it was shown that the asymptotics (\ref{2.10})
 not always can be established in the whole region (\ref{2.11}).
 In some cases the region (\ref{2.11}) should be reduced to the region
  $$
\a:=\frac{x}{n}\in (\a_-,\a_+)\setminus [\b_-,\b_+], ~~~x\in \Z,
 $$
where $\b_-$ is the minimal solution of the equation
$$
  \l(\b)=\l_+,~~~\mbox{where}~~~\l_+:=\sup\{\l:~{\bf E}e^{\l\t}<\infty\},
$$
and $\b_+$ is the maximal solution of this equation.
In the case $\l_+<D(0)$ the both solutions always exist and form the interval
$[\b_-,\b_+]$ of positive length. If $\l_+=D(0)$ and if
$\a=0\in (\a_-,\a_+)$,
then the interval $[\b_-,\b_+]$ degenerate into the point $\b_-=\b_+=0$.
When $\l_+>D(0)$ the interval $[\b_-,\b_+]$ is empty.

Now we can formulate the local theorem for the process
$Z(n)$.  Under the above moment and structural conditions
we will assume also that for $x\in \Z$ the normalized deviations
$\a=\frac{x}{n}$
belong to some fixed compact set
$$
  K\subset (\a_-,\a_+)\setminus [\b_-,\b_+].
$$

Denote
$$
  C(\a):=C_H(1,\a),~~~I(\a):=\sum_{m=1}^\infty e^{\l(\a)m}{\bf P}(\t\ge m),
$$
where the positive continuous in the cone $\mathcal{D}$
function $C_H(\h,\a)$ will be defined below by the formula (\ref{2.26}).

{\bf Theorem 2.1.} {\it   Let us fix the compact
\begin{equation}{\label{2.12}}
K\subset (\a_-,\a_+)\setminus [\b_-,\b_+],
 \end{equation}
and let the condition of admissible heterogeneity holds true:
\begin{equation}{\label{2.13}}
 \mathcal{A}_K\subset (\mathcal{A}_1),~~~\mbox{where}~~~\mathcal{A}_K:=\{(\l,\m)=(\l(\a),\m(\a)):~\a\in K\}.
\end{equation}
Let, in addition to (\ref{2.13}), the following holds true
\begin{equation}{\label{2.14}}
 {\bf P}(\t_1\ge n) =o\left(\frac{1}{\sqrt{n}}e^{-nD(0)}\right)~~~\mbox{when}~~~
 n\to \infty,~~~\mbox{if}~~~\a=0\in K.
\end{equation}
Then when $\a:=\frac{x}{n}\in K$, $x\in \Z$
the representation
\begin{equation}{\label{2.15}}
{\bf P}(Z(n)=x)=\psi_1({\l}(\a),{\m}(\a))\frac{C(\a)}
 {\sqrt{n}}e^{-nD(\a)}I(\a)(1+o(1)),
\end{equation}
holds true, where the remainder function
$o(1)=\v_n(x)$
satisfies the relation
$$
  \lim_{n\to \infty}\sup_{x\in \Z,~
  \frac{x}{n}\in K}|\v_n(x)|=0.
$$
}

Theorem 2.1 is consistent with Theorem 1.1 from \cite{I6}
and extends it for arithmetical case (we mentioned before
that in \cite{I6} the integro-local theorems were established for
non-arithmetical case). It was mentioned in \cite{I6} that a form of the sum
$I(\a)$  explains to some extent an essentiality of the presence of
restriction interval $[\b_-,\b_+]$: if $\a\in (\b_-,\b_+)$
(or, which is the same, $\l(\a)>\l_+$) the sum
$I(\a)$ diverges and the asymptotics
${\bf P}(Z(n)=x)$ will be different.
The essentiality of the condition (\ref{2.13}) explains the presence of the
factor $\psi_1(\l(\a),\m(\a))$
in the right hand side
of (\ref{2.15}).

Let us provide another (equivalent) form of the theorem~2.1:

{\bf Theorem 2.1A.} {\it   Let us fix a point
\begin{equation}{\label{2.16}}
\a_0\in (\a_-,\a_+)\setminus [\b_-,\b_+]
 \end{equation}
and supoose that the condition of admissible heterogeneity holds true:
\begin{equation}{\label{2.17}}
 (\l(\a_0),\m(\a_0))\in (\mathcal{A}_1).
\end{equation}
Adding to (\ref{2.17}), let the following condition holds
\begin{equation}{\label{2.18}}
 {\bf P}(\t_1\ge n) =o\left(\frac{1}{\sqrt{n}}e^{-nD(0)}\right)~~~\mbox{when}~~~
 n\to \infty,~~~\mbox{if}~~~\a_0=0.
\end{equation}
Then for any sequence $x=x_n\in \Z$, such that
$$
  \lim_{n\to \infty}\a=\a_0,~~~\mbox{where}~~~ \a:=\frac{x}{n},
$$
the representation (\ref{2.15}) takes place, where the remainder factor
$o(1)=\v_n(x)$ satisfies
\begin{equation}{\label{2.19}}
  \lim_{n\to \infty}|\v_n(x)|=0.
\end{equation}
}

The equivalence of Theorems 2.1 and 2.1A will be established below.

We now use the following statement:

{\bf Lemma~2.1.} {\it The following two equalities
\begin{equation}{\label{2.20}}
 C(a)I(a)=\frac{1}{\sigma\sqrt{2\pi}},
\end{equation}
\begin{equation}{\label{2.21}}
 \sigma^2=\frac{1}{D''(a)},
\end{equation}
hold true, where (remind)
$$
  \sigma^2:=\frac{1}{{\bf E}\t}{\bf E}(\z-a\t)^2.
$$
}

The statement of Lemma~2.1 can be extracted from \cite{I2}-\cite{I.8}.
However for ``narrative autonomy at the end of \S~3 we provide the proof
of Lemma~2.1. From Theorem~2.1A and Lemma~2.1 it follows

{\bf Corollary~2.1.} {\it Let $\a_0=a$. Then the conditions
(\ref{2.16}), (\ref{2.17}), (\ref{2.18}) of Theorem~2.1A hold true
``automatically" and the relations
$$
D(a)=D'(a)=0,~~~({\l}(a),{\m}(a))=(0,0),~~~
\psi_1({\l}(a),{\m}(a))=1.
$$
hold. Thus, for any sequence
$x=x_n\in \Z$, such that
$$
  \lim_{n\to \infty}\a=a,~~~~\mbox{where}~~~\a=\frac{x}{n},
$$
the following relation
\begin{equation}{\label{2.22}}
{\bf P}(Z(n)=x)\sim \frac{1}{\sigma\sqrt{2\pi n}}e^{-nD(\a)}.
\end{equation}
takes place.}

Corollary~2.1 is proved in \S~3.


Make sure now that the statements of the theorems~2.1 and 2.1A are equivalent, that is,
 the theorem~2.1A follows from the theorem~2.1, and vice versa, the theorem~2.1A
implies the theorem ~2.1.  In other words we establish

\DD~the following implications
\begin{equation}{\label{2.23}}
\mbox{Theorem~2.1}~~~~\Longleftrightarrow~~~~\mbox{Theorem~2.1A}.
\end{equation}
$(i)$. Let the theorem~2.1 holds and let for some point
$\a_0$ the conditions of the theorem~2.1A hold true, that is the relations
(\ref{2.16}), (\ref{2.17}), (\ref{2.18}) hold. Choose $\d>0$ sufficienly
small such that for the interval $K:=[\a_0-\d,~\a_0+\d]$ the relations
(\ref{2.12}), (\ref{2.13}) are satisfied, and additionally, in the case
$\a_0\not =  0$ we have $0\not \in K$. Then, for any sequence
$x=x_n\in \Z$, such that
$\a=\frac{x}{n}\to \a_0$ as $n\to \infty$ according the theorem~2.1 the relation
(\ref{2.15}) takes place, where
the remainder $o(1)=\v_n(x)$ satisfies (\ref{2.19}). Thus, the implication
$$
\mbox{Theorem~2.1}~~~~\Longrightarrow~~~~\mbox{Theorem~2.1A}
$$
is established.

$(ii)$. Let now the theorem~2.1A holds, and suppose that the statement of
the theorem~2.1 is wrong. It means that there exists a compact $K$, which satisfies the relations
(\ref{2.12}), (\ref{2.13}), (\ref{2.14}), and a point
$\a_0\in K$, such that for some sequence of natural numbers $n=n_k\to \infty$ and for some
sequence
$x_{(k)}:=x_{n_k}\in \Z$, such that
$$
    \lim_{k\to \infty}\a_{(k)}=\a_0,~~~\mbox{where}~~~\a_{(k)}:=
    \frac{x_{(k)}}{n_k},
$$
the relation
(\ref{2.15}) takes place, where the reminder
$o(1)=\v_{n_k}(x_{n_k})$
satisfies
\begin{equation}{\label{2.24}}
  \limsup_{k\to \infty}|\v_{n_k}(x_{n_k})|>0.
\end{equation}
But the inequality (\ref{2.24}) is impossible, because it contradicts the relation
(\ref{2.19}), which is true according the theorem~2.1A.
 Thus the implication
$$
\mbox{Theorem~2.1A}~~~~\Longrightarrow~~~~\mbox{Theorem~2.1}
$$
is established.  The equivalence (\ref{2.23}) of the theorem~2.1 and the theorem~2.1A is proved.

Note that Theorem~2.1A is somewhat easier to prove than Theorem~2.1.
Thus in \S~3 we will prove the theorem~2.1A.

{\bf 2.3. Local theorem for renewal function}
$$
  H(B):=\sum_{k=0}^\infty {\bf P}(S_k\in B),~~~B\subset \R^2.
$$
In the arithmetical case for the sequence $(t,x)=(t_n,x_n)\in \Z^2$
    we study the asymptotics
$$
  H(\{t\}\times\{x\}),~~~n\to \infty,
$$
where  $(\h,\a):=\frac{1}{n}(t,x)$ are from some fixed compact set $K$
embedded in the region
$\mathcal{D}$ of analyticity of the function $D(\h,\a)$.
When $\left(\h,\a\right)\in \mathcal{D}$ we set for shortness
(see (\ref{2.10}))
$\left(\widehat{\h},\widehat{\a}\right):=A'\left(\widehat{\l},\widehat{\m}\right)$, where
 \begin{equation}{\label{2.25}}
\left(\widehat{\l},\widehat{\m}\right):=\left(\l\left(\frac{\a}{\h}\right),
 \m\left(\frac{\a}{\h}\right)\right)=
D'\left(\h,\a\right)=\left(-A\left(\m\left(\frac{\a}{\h}\right)\right),
\m\left(\frac{\a}{\h}\right)\right),
\end{equation}
such that the vectors
$\left(\widehat{\l},\widehat{\m}\right)$, $\left(\widehat{\h},\widehat{\a}\right)$
are functions of variable $\frac{\a}{\h}$.
Indeed,
$$
\begin{aligned}
 (\widehat{\l},\widehat{\m})=&\left(\l\left(\frac{\a}{\h}\right), \m\left(\frac{\a}{\h}\right)\right),
\\
 (\widehat{\h},\widehat{\a}) =&
  \left(A'_{\left(1\right)}\left(\widehat{\l},\widehat{\m}\right),
  A'_{\left(2\right)}\left(\widehat{\l},\widehat{\m}\right)\right)
\\
=& \left(A'_{\left(1\right)}\left(\l\left(\frac{\a}{\h}\right),\m\left(\frac{\a}{\h}\right)\right),
 A'_{\left(2\right)}\left(\l\left(\frac{\a}{\h}\right),\m\left(\frac{\a}{\h}\right)\right)\right).
\end{aligned}
$$
It was established in \cite{I6}, \cite{I7} (see
Lemma 3.1 in \S~3 and
Lemma 3.1 in \cite{I6}, \cite{I7}),
 for $\left(\h,\a\right)\in \mathcal{D}$
the minimum over $r\in (0,\infty)$ of the function
$$
  L\left(r\right)=L_{\h,\a}\left(r\right):=r\Lambda\left(\frac{\h}{r},\frac{\a}{r}\right)
$$
is attained at unique point
$$
  r_{\h,\a}=\frac{\h}{\widehat{\h}}=\frac{\a}{\widehat{\a}}.
$$
For $(\h,\a)\in \mathcal{D}$ denote
\begin{equation}{\label{2.26}}
C_H(\h,\a):=\sqrt{\frac{{r}_{\h,\a}}{2\pi} \frac{|\widehat{\Lambda}''|}
{(\h,\a)\widehat{\Lambda}''(\h,\a)^\top}},
\end{equation}
where
$
\widehat{\Lambda}''=\|\widehat{\Lambda}''_{i,j}\|:=
\Lambda''(\h,\a)|_{(\h,\a)=(\widehat{\h},\widehat{\a})},
$
such that
$$
 (\h,\a)\widehat{\Lambda}''(\h,\a)^\top:=
 \h^2\widehat{\Lambda}''_{1,1}+
 2\h\a\widehat{\Lambda}''_{1,2}+
 \a^2 \widehat{\Lambda}''_{2,2}.
$$

%

{\bf Theorem 2.2.} {\it For a fixed compact set
$$
K\subset\mathcal{D}
$$
assume that the condition of admissible heterogeneity holds true:
$$
\mathcal{A}_K\subset  (\mathcal{A}_1),~~~\mbox{where}~~~
\mathcal{A}_K:=\left\{(\l,\m)=\left(\l\left(\frac{\a}{\h}\right),\m\left(\frac{\a}{\h}\right)\right):~(\h,\a)\in K\right\}.
$$
Then for $(t,x)\in \Z^2$,  $(\h,\a):=\frac{1}{n}(t,x)$ the following representation takes place
\begin{equation}{\label{2.27}}
H(\{t\}\times\{x\}) =\frac{1}{\sqrt{n}}\psi_1\left(\l\left(\frac{\a}{\h}\right),~\m\left(\frac{\a}{\h}\right)\right)
C_H(\h,\a)e^{-nD(\h,\a)}(1+o(1)),
\end{equation}
where the reminder $o(1)=\v_n(t,x)$ satisfies
$$
  \lim_{n\to \infty}\sup_{(t,x)\in \Z^2,~\frac{1}{n}(t,x)\in K}|\v_n(t,x)|=0.
$$
}

Theorem 2.2 complements the theorem 3.1 from \cite{I6}, \cite{I7},
where integro-local theorem for renewal function was established
in non-arithmetic case. Theorem 2.2 generalizes the theorem 5 from \cite{I11},
where the only homogeneous case was considered.

Give now the Theorem~2.2A, which statement is equivalent to the statement of
Theorem~2.2, but which is more convenient for proof.

{\bf Theorem 2.2A.} {\it Let for a fixed point
$$
(\h_0,\a_0)\in \mathcal{D}
$$
the condition of admissible heterogeneity holds true:
\begin{equation}{\label{2.28}}
\left(\l\left(\frac{\a_0}{\h_0}\right),~\m\left(\frac{\a_0}{\h_0}\right)\right)\in (\mathcal{A}_1).
\end{equation}
Then for any sequence $(t,x) = (t_n,x_n) \in \Z^2$ such that
$(\h,\a):=\frac{1}{n}(t,x)$ it holds
$$
  \lim_{n\to \infty}(\h,\a)=(\h_0,\a_0),
$$
and the representation (\ref{2.27}) takes place,
where the reminder $o(1)=\v_n(t,x)$ satsfies
$$
  \lim_{n\to \infty}|\v_n(t,x)|=0.
$$
}

The proof of equivalence
\begin{equation}{\label{2.29}}
\mbox{Theorem~2.2}~~~~\Longleftrightarrow~~~~\mbox{Theorem~2.2A}
\end{equation}
repeats the proof of (\ref{2.23}) in the proof of equivalence of
Theorem~2.1 and 2.1A; so we omit the proof of
(\ref{2.29}).
The proof of Theorem~2.2A will be done in \S~3.

\setcounter{section}{3}
\setcounter{equation}{0}

                       \begin{center}
                       {\bf \S~3.
                       The proof of the main statements.   }
                       \end{center}

{\bf 3.1. \DD~of Theorem~2.2A.}
The initial statement for the proof of Theorem 2.2A is
the local theorem for the sum $S_n$ in {\it arithmetical homogeneous case}.
Theorem~2.3.2, see \cite{I1}, p. 72, implies the following version of
local theorem for the sums $S_n$ of random vectors in arithmetical case:

{\bf Theorem 3.1.} {\it  Let us fix a compact
$K\subset \mathcal{L}$,  a point $(t,x)\in \Z^2$
such that $(\g,\b):=(\frac{t}{n},\frac{x}{n})\in K$.
Then with
 $$
 C_1(\g,\b):=\frac{\sqrt{|\Lambda''(\g,\b)|}}{2\pi},
$$
it holds
 \begin{equation}{\label{3.1}}
{\bf P}(T_n=t,~Z_n=x)=
\frac{1}{n }
C_1(\g,\b)
e^{-n\Lambda(\g,\b)}(1+o(1)),
 \end{equation}
as $n\to \infty$, where remainder $o(1)=\v_n(t,x)$ is uniform on
   $(\g,\b)\in K$:
$$
  \lim_{n\to \infty} \sup_{(t,x)\in \Z^2,~(\frac{t}{n},\frac{x}{n})\in K}
  |\v_n(t,x)|=0.
$$
   }

The asumptotics of renewal measure  $H(B)$ will be find below using (\ref{3.1})
for receding sets
 $$
 B=B_{t,x}:=\{t\}\times\{x\},
 $$
that is, the statement (\ref{2.27}) of Theorem~2.2A will be established.

  For the proof of the theorem~2.2A we need

{\bf Lemma 3.1.}(\cite{I6}, \cite{I7}) {\it
Let $(\h,\a)\in {\mathcal{D}}$, and the vectors
$(\widehat{\l},\widehat{\m})$, $(\widehat{\h},\widehat{\a})$
are defined in (\ref{2.25}). Then:

$I$. The minimum of the function $L(r)=L_{\h,\a}(r):=r\Lambda(\frac{\h}{r},\frac{\a}{r})$
on the set $r > 0$ is attained at the unique point
$$
  r_{\h,\a}=\frac{\h}{\widehat{\h}}=
 \frac{\a}{\widehat{\a}}.
$$
For the functions $\widehat{\l}$, $\widehat{\m}$ the following representation takes place
$$
\widehat{\l}=\l\left(\frac{\h}{r_{\h,\a}},\frac{\a}{r_{\h,\a}}\right),~~~
\widehat{\m}=\m\left(\frac{\h}{r_{\h,\a}},\frac{\a}{r_{\h,\a}}\right),
$$
where the functions $\l(\cdot,\cdot)$, $\m(\cdot,\cdot)$
are the solutions of the system (\ref{2.2}).
Herewith
\begin{equation}{\label{3.2}}
  L'(r_{\h,\a})=0,~~~L''(r_{\h,\a})=  \frac{1}{r_{\h,\a}}(\widehat{\h},\widehat{\a})
\widehat{\Lambda}''(\widehat{\h},\widehat{\a})^\top>0,
\end{equation}
where
$$
\widehat{\Lambda}'':=\Lambda''(\h,\a)|_{(\h,\a)=(\widehat{\h},\widehat{\a})}.
$$

$II$. For a given  $\v>0$ there exists $\v_1>0$ such that
$$
\min_{|r-r_{\h,\a}|\ge \v}L(r)\ge L(r_{\h,\a})+\v_1.
$$
}

\DD~of Theorem~2.2A (based on Theorem~3.1 and Lemma~3.1) we divide into two steps.

$I$. Denote $H_0(B)$ the renewal measure for homogeneous case, when the vectors
$\x_1=(\t_1,\z_1)$, $\x=(\t,\z)$ (i.e. all terms
$\x_i=(\t_i,\z_i)$ when $i\ge 1$) have the same distribution. In the first step
we prove the statement of Theorem~2.2A for homogeneous case.
Here it is sufficient to prove that for any fixed
$(\h_0,\a_0)\in \mathcal{D}$ and for any sequence
$(t,x)\in \Z^2$, such that
 $(\h,\a):=\frac{1}{n}(t,x)\to (\h_0,\a_0)$
as $n\to \infty$ the equality
\begin{equation}{\label{3.3}}
H_0(\{t\}\times\{x\}) =\frac{1}{\sqrt{n}} C_H(\h,\a)e^{-nD(\h,\a)}(1+o(1)),
\end{equation}
takes place. Let us prove now (\ref{3.3}).
We have
\begin{equation}{\label{3.4}}
H_0(B_{t,x})=\sum_{k=1}^\infty {\bf P}(T_k=t,~Z_k=x).
 \end{equation}
The series in (\ref{3.4}) we split into three sums
$$
  \sum_{k\in \mathcal{K}_1},~~~\sum_{k\in \mathcal{K}_2},~~~\sum_{k\in \mathcal{K}_3},
$$
over regions
$$
   \mathcal{K}_1:=\{1\le k<c_1n\},~~~ \mathcal{K}_2:=\{c_1n\le k\le c_2n\},~~~
   \mathcal{K}_3:=
  \{ k> c_2n\},
$$
where $c_1<c_2$ are chosen such that
$$
  r_{\h_0,\a_0}\in (c_1,c_2),
$$
and the difference $c_2-c_1$ is so small that for some compact
$K\subset \mathcal{L}$ and for all sufficiently large
$n$ when
$k\in \mathcal{K}_2$
it holds that
  $$
   (\g,\b):=\left(\frac{t}{k},\frac{x}{k}\right)\in K.
 $$
It is obvious that it is always possible.
Thus, according Theorem~3.1 for $r:=\frac{k}{n}\in [c_1,c_2]$
   $$
   \sum_{k\in \mathcal{K}_2}= \sum_{k\in \mathcal{K}_2}
   \frac{C_1(\g,\b)}{nr}
   e^{-k\Lambda(\g,\b)}(1+o(1))=
 $$
 $$
   \sum_{r\in [c_1,c_2]}\frac{1}{n}
  \frac{C_1(\frac{\h}{r},\frac{\a}{r})}{r}
   e^{-nr\Lambda(\frac{\h}{r},\frac{\a}{r})}(1+o(1)).
 $$
Since the variation of the function $r\Lambda(\frac{\h}{r},\frac{\a}{r})=L(r)$ is $o(\frac{1}{n})$
on the neighborhood of the point $r=r_{\h,\a}\in [c_1,c_2]$ on the interval of length
$\frac{1}{n}$, then
$$
   \sum_{k\in \mathcal{K}_2}=\int_{c_1}^{c_2}\frac{1}{r}
      C_1\left(\frac{\h}{r},\frac{\a}{r}\right)e^{-nL(r)}dr (1+o(1)).
$$
By virtue of the equality $L(r_{\h,\a})=D(\h,\a)$ and the famous Laplace method
of calculating such integrals we obtain
(for  $(\widehat{\h},\widehat{\a})=
(\frac{\h}{r_{\h,\a}},\frac{\a}{r_{\h,\a}})$)
$$
  \sum_{k\in \mathcal{K}_2}=
  \frac{C_1(\widehat{\h},\widehat{\a})\sqrt{2\pi}}
  {\sqrt{nL''(r_{\h,\a})}r_{\h,\a}}e^{-nD(\h,\a)}(1+o(1)).
$$
According Lemma~3.1 (see the right formula in (\ref{3.2})) the right-hand side of this
equality coincides with the right-hand side of
(\ref{3.3}) (that is, with the right part of (\ref{2.27}) when $\psi_1=\psi$).

Estimate now $\sum_{k\in \mathcal{K}_1}$ and
$\sum_{k\in \mathcal{K}_3}$. By virtue of multidimensional exponential
inequality of Chebyshev's type (see \cite{I10}, or Theorem~1.3.2 in \cite{I1})
$$
  {\bf P}(T_k=t, ~Z_k=x) \le
   \exp\skk{-nr\Lambda\left(\frac{\h}{r},\frac{\a}{r}\right)},~~~r=\frac{k}{n}.
$$
According the statement $II$ of Lemma~3.1 there exist $\g>0$ and $n_0<\infty$
such that for all $n\ge n_0$, $k\in \mathcal{K}_1\cup\mathcal{K}_3$
it holds that
$$
  r\Lambda\left(\frac{\h}{r},\frac{\a}{r}\right)
  \ge D(\h,\a)+\g.
$$
Thus, when $n\ge n_0$ we have
\begin{equation}{\label{3.5}}
{\bf P}(T_k=t, ~Z_k=x) \le e^{-n(D(\h,\a)+\g)}.
 \end{equation}
Therefore, by (\ref{3.5})
$$
\sum_{k\in \mathcal{K}_1}\le c_1n e^{-n(D(\h,\a)+\g)} =
o\left( \frac{1}{\sqrt{n}}e^{-nD(\h,\a)}\right).
$$
In order to estimate the sum $\sum_{k\in \mathcal{K}_3}$
we split it into two parts:
the sum $\sum'$ over the region $c_2n\le k\le n^2$
and the sum $\sum''$ over the region $ k > n^2$.
The sum $\sum'$ can be estimated in the same way as the sum
$\sum_{k\in \mathcal{K}_1}$.
For the sum   $\sum''$ we have
$$
  {\sum}'' \le \sum_{k > n^2}{\bf P}(T_k \le 2t)\le
  \sum_{k > n^2}e^{-k\Lambda_\t(\frac{t}{k})},
$$
where for $k>n^2\to \infty$,
$$
  \Lambda_\t\left(\frac{t}{k}\right)=\Lambda_\t\left(\frac{\h n}{k}\right)\ge
  \Lambda_\t\left(\frac{\h}{n}\right)\to \infty,
$$
because
$\Lambda_\t(0)=\infty$. Thus
$$
{\sum}'' \le \sum_{k\ge n^2} e^{-k\Lambda_\t(\frac{\h}{n})}\le
 2 e^{-n^2\Lambda_\t(\frac{\h}{n})}=
 o\left(\frac{1}{\sqrt{n}}e^{-nD(\h,\a)}\right).
$$
It proves the theorem~2.2A in homogeneous case
(that is the relation (\ref{3.3})).

$II$. On the second step of the proof we represent the
renewal measure in the following way
$$
  H(B_{t,x})=H(\{t\}\times \{x\})=
  {\bf E}H_0(\{t-\t_1\}\times  \{x-\z_1\}),
$$
and we make use of the following statement (compare with Lemma~4.2 (\cite{I7}),
which was established in the proof of Theorem~3.1 in
 \cite{I7}):

{\bf Lemma 3.2.}
{\it
Let conditions of Theorem~2.2A hold true. Then
for some $c>0$,
$C<\infty$, $n_0<\infty$ for all $n\ge n_0$ it holds that}
\begin{equation}{\label{3.6}}
|H(\{t\}\times\{x\})-
{\bf E}(H_0(\{t-\t_1\}\times
 \{x-\z_1\});~~
  |\x_1|\le \ln^2n)|
  \le Ce^{-nD(\h,\a)-c\ln^2n}.
 \end{equation}

\DD. Let for $(t,x)\in \Z^2$, $t\ge 1$
$$
   B_{t,x}=\{t\}\times \{x\}.
$$
Since $(0,0)\not \in B_{t,x}$, then
$$
{\bf P}((T_0,Z_0)\in B_{t,x})=0,
$$
thus, we have
$$
 H(B_{t,x})=\sum_{k=1}^\infty {\bf P}((T_k,~Z_k)\in B_{t,x})=
$$
$$
  \sum_{k=0}^\infty {\bf P}((\t_1+T'_k,~\z_1+Z'_k)\in B_{t,x}),
$$
where $(T'_k,Z'_k)$, $k=0,1,\cdots$ is the independent on $(\t_1,\z_1)$ sequence
of sums of homogeneous independent terms. Thus the renewal measure $B_{t,x}$
in non-homogeneous case can be represented as
\begin{equation}{\label{3.7}}
  H(B_{t,x})=H(\{t\}\times \{x\})={\bf E}H_0(\{t-\t_1\}\times
  \{x-\z_1\}).
\end{equation}
According (\ref{3.7}) we have
$$
H(\{t\}\times \{x\})=
{\bf E}(H_0(\{t-\t_1\}\times
  \{x-\z_1\});~|\x_1|\le \ln^2n)+
$$
$$
{\bf E}(H_0(\{t-\t_1\}\times  \{x-\z_1\});~|\x_1|> \ln^2n).
$$
 Thus
$$
\Bigl|H(\{t\}\times \{x\})-{\bf E}(H_0(\{t-\t_1\}\times  \{x-\z_1\});~~
|\x_1|\le \ln^2n)\Bigr|=
  \sum_{k=1}^\infty P_k,
$$
{ where}
 $$
  P_k:={\bf P}(T_k = t,~Z_k = x,~|\x_1|>\ln^2n),
  $$
  and for the proof of Lemma~3.2 it is sufficient to prove that
\begin{equation}{\label{3.8}}
\sum_{k=1}^\infty P_k\le Ce^{-n{\bf D}(\h,\a)-c\ln^2n}.
 \end{equation}
For all sufficiently large $n$ it holds that
$(\widehat{\l},\widehat{\m})\in (\mathcal{A})$. Thus for such
$n$ for any $k\ge 1$     
we have
  $$
  P_k={\bf E}(e^{-\widehat{\l}T_k-\widehat{\m}Z_k
  +\widehat{\l}T_k+\widehat{\m}Z_k};~T_k = t,~Z_k = x,~|\x_1|>\ln^2n)
  \le
  $$
 $$
  \exp \{-\widehat{\l}t-\widehat{\m}x\}
  {\bf E}(e^{\widehat{\l}T_k+\widehat{\m}Z_k};|\x_1|>\ln^2n)=
  $$
  $$
  \exp \{-n{ D}(\h,\a)\}
  {\bf E}(e^{\widehat{\l}\t_1+\widehat{\m}\z_1};~|\x_1|>\ln^2n)
  \prod_{j=2}^k {\bf E}e^{\widehat{\l}\t_j+\widehat{\m}\z_j}\le
 $$
 $$
  \exp \{-n{ D}(\h,\a)\}
  {\bf E}(e^{\widehat{\l}\t_1+\widehat{\m}\z_1};~|\x_1|>\ln^2n).
 $$
 {  In the last inequality we took advantage of the
 fact that for all sufficiently large $n$ the vector $(\widehat{\l},\widehat{\m})$
belongs to the boundary of the set
 $\mathcal{A}^{\le 0}$, and for $j\ge 2$ it holds that
 $$
  {\bf E}e^{\widehat{\l}\t_j+\widehat{\m}\z_j} = 1.
 $$

 Applying further the H\"older inequality we obtain that for
 $p>0$, $q>0$, $p+q=1$}
 $$
  {\bf E}(e^{\widehat{\l}\t_1+\widehat{\m}\z_1};~|\x_1|>\ln^2n)\le
  \Bigl({\bf E}e^{\frac{1}{p}\widehat{\l}\t_1+\frac{1}{p}\widehat{\m}\z_1}
  \Bigr)^p
  {\bf P}^q\left(|\x_1|>\ln^2n\right).
  $$
 { Since by virtue of the condition of admissible heterogeneity  (\ref{2.28})
  for all sufficiently large $n$  it holds
 $(\widehat{\l},\widehat{\m})\in (\mathcal{A}_1)$. Then, choosing parameter
$\frac{1}{p}>1$ close enough to $1$, we obtain
 $$
    {\bf E} e^{
    \frac{1}{p}\widehat{\l}\t_1+\frac{1}{p}\widehat{\m}\z_1}
    \le C_0<\infty.
 $$
 For  $q=1-p$ and for some $C_1<\infty$, $c_1>0$ we have
 (by virtue of $[{\bf C}_0]$)}
$$
  {\bf P}^q(|\x_1|> \ln^2n)\le C_1e^{-c_1\ln^2n}.
$$
 { Thus for any $k\ge 1$ the following inequality holds true
 \begin{equation}{\label{3.9}}
P_k\le
C_0^pC_1e^{-n{ D}(\h,\a)-c_1\ln^2n}.
 \end{equation}
Further,
\begin{equation}{\label{3.10}}
\sum_{k=1}^\infty P_k\le n^2 \sup_{k\ge 1}P_k+
\sum_{k\ge n^2}P_k,
 \end{equation}
where the first term in right-hand side can be estimated using
(\ref{3.9}), and the second one using inequalities
\begin{equation}{\label{3.11}}
  \limsup_{n\to \infty}\frac{1}{n}\ln \sum_{k\ge n^2}P_k\le
  \limsup_{n\to \infty}\frac{1}{n}\ln \sum_{k\ge n^2}
{\bf P}(T_k=t)=-\infty.
 \end{equation}
From (\ref{3.9})---(\ref{3.11}) it follows that
for some $c>0$, $C<\infty$ the inequality (\ref{3.8}) holds true.
Lemma~3.2 is proved.

 By (\ref{3.6})
 for the proof of (\ref{2.27}) in general case it is sufficient to find an asymptotics of the
 mean
 $$
   {\bf E}(H_0(\{t-\t_1\}\times  \{x-\z_1\});~~
|\x_1|\le \ln^2n).
 $$
 According the result of the first step and by (\ref{3.6}) we have
 $$
H(\{t\}\times \{x\})=
$$
\begin{equation}{\label{3.12}}
 \frac{1}{\sqrt{n}}{\bf E}\left(C_1\left(\h-\frac{\t_1}{n},\a-\frac{\z_1}{n}\right)
e^{-nD\left(\h-\frac{\t_1}{n},\a-\frac{\z_1}{n}\right)}(1+o(1));~~|\x_1|\le \ln^2n\right)+
 \end{equation}
$$
o\left(\frac{1}{\sqrt{n}}e^{-nD(\h,\a)}\right).
$$
Considering the fact that in the region $|\x_1|\le \ln^2n$, when $n\to \infty$,
it holds that
$$
  -nD\left(\h-\frac{\t_1}{n},\a-\frac{\z_1}{n}\right)=-nD(\h,\a)+
  \widehat{\l}\t_1+\widehat{\m}\z_1+o(1),
$$
we obtain that the right-hand side of (\ref{3.12}) coincides with right-hand side of (\ref{2.27}).
 Theorem~2.2A is proved.

{\bf 3.2. \DD~of Theorem~2.1A.}
We have
$$
  P_n:={\bf P}(Z(n)=x)=
  $$
  $$
  {\bf P}(Z_0=x,~\n(n)=0)+\sum_{k=1}^n
  {\bf P}(Z_k=x,~\n(n)=k)=
  {\bf I}_{\{x=0\}}{\bf P}(\t_1\ge n)+R_n,
$$
where
$$
 K_n:=\sum_{k=1}^n
  {\bf P}(Z_k=x,~\n(n)=k)=
  \sum_{k=1}^n
  {\bf P}(Z_k=x,~T_n<n,~T_{k+1}\ge n).
$$
Using the renewal function
$H(\{\cdot\}\times
\{x\})$ the term $R_n$ we represent as
$$
 K_n=\sum_{m=0}^{n-1}\sum_{k=1}^n
  {\bf P}(Z_k=x,~T_k=m,~T_{k+1}\ge n)=
  $$
  $$
  \sum_{m=0}^{n-1}\sum_{k=1}^n
  {\bf P}(Z_k=x,~T_k=m){\bf P}(\t\ge n-m)=
$$
$$
 \sum_{m=0}^{n-1}H(\{m\}\times\{x\}){\bf P}(\t\ge n-m).
$$
Thus, we obtain the representation
\begin{equation}{\label{3.13}}
P_n={\bf I}_{\{x=0\}}{\bf P}(\t_1\ge n)+
\sum_{m=0}^{n-1}Q_n(m),
\end{equation}
where
$$
Q_n(m):=H(\{m\}\times\{x\}){\bf P}(\t\ge n-m),~~~m\in \{0,\cdots,n-1\}.
$$

We use the following statement (compare with Lemma~3.2):

{\bf Lemma~3.3.} {\it Let the conditions of Theorem~2.1A hold true. Then
for some $c>0$, $C<\infty$, $n_0<\infty$, for all $n\ge n_0$
the following inequality holds
$$
R_n:={\bf P}(Z(n)=x, \t_{\n(n)+1}\ge \ln^2n)\le Cne^{-nD(\a)-c\ln^2n}.
$$
}
We will prove Lemma~3.3 later in this section, but now, using this lemma, we finish the proof
of Theorem~2.1A. Find the asymptotics
$$
  L_n:=\sum_{k=1}^{[\ln^2n ]}H(\{n-k\}\times\{x\}){\bf P}(\t\ge k).
$$
By Theorem~2.2A, we obtain
$$
  L_n=\sum_{k=1}^{[\ln^2n ]}\frac{1}{\sqrt{n}} \psi_1(\l(\a),\m(\a))C_H(1,\a)
  e^{-nD\left(1-\frac{k}{n},\a\right)}{\bf P}(\t\ge k)(1+o(1)).
$$
Since in the interval $1\le k\le [\ln^2n ]$ it holds
$$
 -nD\left(1-\frac{k}{n},\a\right)=-D(1,\a)+(\l(\a)+o(1))k,
$$
thus, from the last we obtain
$$
  L_n=\frac{1}{\sqrt{n}} \psi_1(\l(\a),\m(\a))C_H(1,\a)
  e^{-nD(1,\a)}\sum_{k=1}^{[\ln^2n ]}e^{\l(\a)k}{\bf P}(\t\ge k)(1+o(1))=
$$
$$
  \frac{1}{\sqrt{n}} \psi_1(\l(\a),\m(\a))C_H(1,\a)
  e^{-nD(1,\a)}I(\a)(1+o(1)).
$$
In this way we obtained the asymptotic which coincides with right-hand side of
the relation (\ref{2.15}). Since from (\ref{3.13}) and from Lemma~3.3 it follows
that for all sufficiently large $n$
$$
  |P_n-L_n|\le R_n\le Cne^{-nD(\a)-c\ln^2n},
$$
then we obtain the required relation
(\ref{2.15}).


It remains to provide

\DD~Lemma~3.3. First we estimate each term in the right-hand side of the equality
$$
 {\bf P}(Z(n)=x,~~~\t_{\n(n)+1}\ge \ln^2n)=
 $$
\begin{equation}{\label{3.14}}
 {\bf I}_{\{x=0\}}
 {\bf P}(\t_1\ge n)+
 \sum_{k=1}^n {\bf P}(Z_k=x,~~T_k<n,~T_{k+1}\ge n,~~B_n) =:\sum_{k=0}^nR_k(n),
\end{equation}
where $B_n:=\{\t_{k+1}\ge \ln^2n\}$.
According the conditions of Theorem~2.1A we have for some $c>0$, $C<\infty$
$$
  R_0(n)={\bf I}_{\{x=0\}}{\bf P}(\t_1\ge n)\le Ce^{-nD(\a)-cn}.
$$
For $k\ge 1$ further we have
$$
R_k(n):={\bf P}(Z_k=x,~~T_k<n,~T_{k+1}\ge n,~~B_n)=
$$
$$
e^{-n(\l(\a)+\m(\a)\frac{x}{n})}
{\bf E}(e^{n\l(\a)+\m(\a)Z_k};~Z_k=x,~~T_k<n,~T_{k+1}\ge n,~~B_n).
$$
Applying an absolutely-continous mapping which transforms the distribution of the vectors
$(\t_1,\z_1)$, $(\t,\z)$ in the distribution
$$
  {\bf P}((\widehat{\t}_1,\widehat{\z}_1)\in B):=
  \frac{1}{\psi_1(\l(\a),\m(\a))}{\bf E}(e^{\l(\a)\t_1+\m(\a)\z_1};~
  (\t_1,\z_1)\in B),
  $$
  $$
  {\bf P}((\widehat{\t},\widehat{\z})\in B):=
  {\bf E}(e^{\l(\a)\t+\m(\a)\z};~
  (\t,\z)\in B),
$$
with the natural interpretation of the notations $\widehat{T}_k$,
$\widehat{Z}_k$ we obtain the equality
$$
 R_k=
 e^{-nD(\a)}\psi_1(\l(\a),\m(\a)) E_k,
$$
where
$$
E_k:={\bf E}(e^{
\l(\a)(n-\widehat{T}_k)};~
\widehat{Z}_k=x,~~\widehat{T}_k<n,~
\widehat{T}_k+\t_{k+1}\ge n,~~B_n).
$$
Estimate from above $E_k$.

If $\l(\a)\le 0$, then on the event $\{\widehat{T}_k<n\}$
we have $e^{
\l(\a)(n-\widehat{T}_k)}\le 1$ and
$$
E_k\le {\bf P}(B_n)={\bf P}(\t\ge \ln^2n) \le Ce^{-c\ln^2n}.
$$
If $\l(\a)>0$, then on the event $\{\widehat{T}_k+\t_{k+1}\ge n\}$
we have
$$
  e^{\l(\a)(n-\widehat{T}_k)}=e^{\l(\a)(n-\widehat{T}_k-\t_{k+1})+
  \l(\a)\t_{k+1}}\le e^{\l(\a)\t_{k+1}}
$$
and
\begin{equation}{\label{3.15}}
 E_k\le
 {\bf E}(e^{\l(\a)\t_{k+1}},~~B_n)={\bf E}(e^{\l(\a)\t},~~
 \t\ge \ln^2n).
\end{equation}
By virtue of the H\"older inequality for $p>0$, $q>0$, $p+q=1$, we have
\begin{equation}{\label{3.16}}
 {\bf E}(e^{\l(\a)\t_{k+1}},~~B_n)={\bf E}(e^{\l(\a)\t},~~
 \t\ge \ln^2n)\le {\bf E}^pe^{\frac{\l(\a)}{p}\t}{\bf P}^q(\t\ge \ln^2n).
\end{equation}
Due to the conditions of Theorem~2.1A there exist $\v_0>0$, $p_0\in (0,1)$, $C_0<\infty$
such that for all sufficiently large $n$
$$
  \frac{\l(\a)}{p_0} \le \l_+-\v_0,~~~{\bf E}^pe^{\frac{\l(\a)}{p}\t}\le C_0,
  $$
thus, from (\ref{3.15}), (\ref{3.16}) it follows that for $q_0=1-p_0$
$$
 E_k\le C_0{\bf P}^{q_0}(\t\ge \ln^2n).
$$
We obtained in this case again for some
$c>0$, $C<\infty$ the inequality
$$
 E_k\le   Ce^{-c\ln^2n}.
$$

Thus,  for some
$c>0$, $C<\infty$, $n_0<\infty$ and for all $n\ge n_0$
each term in the right-hand side of (\ref{3.14})
satisfies inequality
\begin{equation}{\label{3.17}}
R_k(n) \le  Ce^{-nD(\a)-c\ln^2n}.
\end{equation}
By (\ref{3.17}) we have
$$
{\bf P}(Z(n)=x,~~~\t_{\n(n)+1}\ge \ln^2n)\le (n+1)C e^{-nD(\a)-c\ln^2n}.
$$
Lemma~3.3 is proved. Theorem~2.1A is proved.

{\bf 3.3. \DD~of Corollary~2.1.} Due to the fact that
for random vector $\x=(\t,\z)$ the condition $[{\bf C}_0]$ holds,
the point $(\l,\m)=(0,0)$ belongs to the region
$\mathcal{A}$ of analiticity of the function $A(\l,\m)$.
Thus, the point $(a_\t,a_\z)=A'(\l,\m)|_{(\l,\m)=(0,0)}$ belongs to the region
$\mathcal{L}$ of analyticity of the deviation function $\Lambda(\h,\a)$.
Hence the point $a:=\frac{a_\z}{a_\t}$ lies in the cone
$\mathcal{D}$ (the region of analyticity of the function $D(\a)$).
Thus due to the lemma~2.1 we have
$$
  (\l(a),\m(a))=(0,0),~~~D(a)=\l(a)+\m(a)a=0,~~~D(a)=\m(a)=0.
$$
Since the random vector $\x_1=(\t_1,\z_1)$ satisfies the condition
$[{\bf C}_0]$, then for $\a_0=a$
the conditions (\ref{2.16}), (\ref{2.17}), (\ref{2.18}) of the theorem~2.1A
hold true ``automatically''.
Thus, applying Theorem~2.1A, for $\a_0=a$ we obtain the statement (\ref{2.22})
of Corollary~2.1. Corollary~2.1 is proved.

{\bf 3.4. \DD~of Lemma~2.1.}
I. Let us prove (\ref{2.21}). Remind, that for the points
$\a$ from some neigborhood of point $\a_0=a$ it holds
$$
D(\a)=\m(\a)\a-A(\m(\a)),
$$
where the function $A(\m)$ is the solution of equation
\begin{equation}{\label{3.18}}
\psi(-A(\m),\m)=1,~~~\psi(\l,\m):={\bf E}e^{\l\t+\m\z},
\end{equation}
and the function $\m(\a)$ is the inverse for the function $A'(\m)$,
that is, it satisties the equation
$$
A'(\m(\a))=\a.
$$
Make sure that
\begin{equation}{\label{3.19}}
A(0)=0,~~~A'(0)=a,~~~A''(0)=\sigma^2,
\end{equation}
\begin{equation}{\label{3.20}}
D(a)=0,~~~D'(a)=0,~~~D''(a)=\frac{1}{\sigma^2}.
\end{equation}

Equalities (\ref{3.19}) are obtained from the equation (\ref{3.18})
by differentiating it the required number of times:

$(0)$ when $\m=0$, due to $\psi(0,0)=1$
we have $A(0)=0$;

$(1)$ we differentiate equation (\ref{3.18}) one time when
$\m=0$, we have
$$
  \psi'_1(0,0)(-A'(0))+\psi'_2(0,0)=0;
$$
thus
$$
  A'(0)=\frac{\psi'_2(0,0)}{\psi'_1(0,0)}=a;
$$

$(2)$ we differentiate equation (\ref{3.18}) two times when
$\m=0$, we have
$$
  \psi''_{1,1}(0,0)(-A'(0))^2-2\psi''_{1,2}(0,0)A'(0)+
  \psi''_{2,2}(0,0)-\psi'_{1}(0,0)A''(0)=0;
$$
thus
$$
  A''(0)=\frac{1}{\psi'_1(0,0)}[{\bf E}\z^2-2a{\bf E}\z\t+a^2{\bf E}\t^2]=
  \sigma^2.
$$
Equalities (\ref{3.19}) are established. Since the functions $D'(\a)$ and $A'(\m)$
are reciprocal, then (\ref{3.19}) implies
(\ref{3.20}).  The equality (\ref{2.21}) is proved.

II. Let us prove (\ref{2.20}).  Remind that
$$
  C(\a):=C_H(1,\a),~~~I(\a):=\int_0^\infty e^{\l(\a)y}{\bf P}(\t>y)dy,~~~
  \l(\a)=-A(\m(\a)),
$$
where
$$
C_H(\h,\a)=
\sqrt{
\frac{
\hat{\h}|\Lambda''(\hat{\h},\hat{\a})|
}
{
2\pi\h Q(\hat{\h},\hat{\a})
}
},
$$

$$
Q(\hat{\h},\hat{\a})=(\hat{\h},\hat{\a})
\Lambda''(\hat{\h},\hat{\a})(\hat{\h},\hat{\a})^T,
$$
\begin{equation}{\label{3.21}}
 \left(\hat{\h},\hat{\a}\right)=A'\left(\l\left(\frac{\a}{\h}\right),\m\left(\frac{\a}{\h}\right)\right).
\end{equation}
Make sure that the following equalities hold true
\begin{equation}{\label{3.22}}
 I(a)={\bf E}\t,~~~C(a)=\frac{1}{\sigma{\bf E}\t\sqrt{2\pi}}.
\end{equation}

The first equality in (\ref{3.22}) is obvious.  From the definition
(\ref{3.21}) it follows that, when $(\h,\a)=(1,a)$:
$$
 (\hat{\h},\hat{\a})=A'(\l(a),\m(a))=A'(0,0)=(a_\t,a_\z).
$$
Therefore, (see the second equality in \cite{I1})
$$
 \Lambda''(\hat{\h},\hat{\a})=\Lambda''(a_\t,a_\z)= M^{-1},
$$
where
$$
 M=\|M_{i,j}\|_{i,j=1,2}
$$
is the covariance matrix of random vector $(\t,\z)$:
$$
  M_{1,1}={\bf E}\t_0^2,~~~M_{1,2}=M_{2,1}={\bf E}\t_0\z_0,~~~
  M_{2,2}={\bf E}\z_0^2,
$$
where $\t_0:=\t-a_\t$, $\z_0:=\z-a_\z$.
Thus,
$$
\Lambda''(\hat{\h},\hat{\a})=\frac{1}{|M|}C,~~~C=\|C_{i,j}\|_{i,j=1,2},
$$
where $|M|$ is determinant of the matrix $M$,
$$
  C_{1,1}=M_{2,2},~~~
  C_{2,2}=M_{1,1},~~~C_{1,2}=C_{2,1}=-M_{1,2}.
$$
Next, we have
$$
 |\Lambda''(\hat{\h},\hat{\a})|=\frac{1}{|M|},
$$
$$
  Q(\hat{\h},\hat{\a})=\frac{1}{|M|}(a_\t)^3 {\bf E}(\z_0-a\t_0)^2=
  \frac{1}{|M|}(a_\t)^3 {\bf E}(\z-a\t)^2=
  \frac{1}{|M|}(a_\t)^3\sigma^2.
$$
Thus,
$$
  C(a)=\sqrt{\frac{a_\t\frac{1}{|M|}}{\frac{2\pi (a_\t)^3 \sigma^2}{|M|}}}=
  \frac{1}{a_\t \sigma\sqrt{2\pi}}.
$$
The second equality in (\ref{3.22}) is established.
Lemma~2.1 is proved.

The author thanks E.I.Prokopenko for proof of Lemma~2.1 and
A.V.Logachov for valuable comments.


\begin{thebibliography}{99}

\bibitem{I1}
A.A. Borovkov,
{\it Asymptotic analysis of random walks. Rapidly decreasing distributions of increments,}  Moscow:Fizmatlit, 2013.




\bibitem{I2}
A.A. Borovkov,  A.A. Mogulskii,
{ \it Large deviation principles for the finite-dimensional distributions of compound renewal processes}, Sib. Math. J., {\bf 56}:1 (2015),  28--53.


\bibitem{I3}
A.A. Borovkov,  A.A. Mogulskii,
Large Deviation Principles for Trajectories of Compound Renewal Processes. I,
Theory Probab. Appl., {\bf 60}:2 (2016), 207--224.

 \bibitem{I4}
A.A. Borovkov,  A.A. Mogulskii,
Large Deviation Principles for Trajectories of Compound Renewal Processes. II,
 heory Probab. Appl., {\bf 60}:3 (2016), 349--366.




 \bibitem{I5}
A.A. Borovkov,
Large deviation principles in boundary problems for compound renewal processes,  Siberian Mathematical Journal, {\bf 57}:3 (2016), 562--595




 \bibitem{I6}
 A.A. Borovkov,  A.A. Mogulskii,
{\it Integro-local limit theorems for compound renewal processes with Cramer's condition I}, Siberian Mathematical Journal, {\bf 59}:3 (2018), 491--514.

 \bibitem{I7}
 A.A. Borovkov,  A.A. Mogulskii,
{\it Integro-local limit theorems for compound renewal processes with Cramer's condition II}, Siberian Mathematical Journal, {\bf 59}:4 (2018), 731--750.


\bibitem{I.8}
A.A. Mogulskii, E.I. Prokopenko,
{\it Integro-local theorems for multidimensional compound renewal processes, when Cramer's condition holds. I,}
Siberian Electronic Mathematical
Reports, {\bf 15} (2018), 475--502.

\bibitem{I.9}
A.A. Mogulskii, E.I. Prokopenko,
{\it Integro-local theorems for multidimensional compound renewal processes, when Cramer's condition holds. III,} Siberian Electronic Mathematical
Reports, {\bf 15} (2018), 475--502.


\bibitem{I.10}
A.A. Mogulskii, E.I. Prokopenko,
{\it Integro-local theorems for multidimensional compound renewal processes, when Cramer's condition holds. III,} Siberian Electronic Mathematical
Reports, { \bf 15} (2018), 528--553.

\bibitem{I8}
D.P. Cox, W. L. Smith,
{\it Renewal Theory [Russian translation]},  Moscow:Soy. Radio, 1967.

\bibitem{I9}
S. Asmussen, H. Albrecher,
{\it  Ruin Probabilities}, Second Edition,
Word Scientifics,  2010.



\bibitem{I10}
A.A. Borovkov,  A.A. Mogulskii,
Chebyshev-Type Exponential Inequalities for Sums of Random Vectors and for Trajectories of Random Walks,
Theory Probab. Appl., { \bf 56}:1 (2012), 21--43.


  \bibitem{I11}
A.A. Borovkov,  A.A. Mogulskii,
 {\it The second rate function and the asymptotic problems of renewal and hitting the boundary for multidimensional random walks},
Sib. Math. J., {\bf 37}:4, (1996), 647--682.







%
%
%



 \end{thebibliography}
 \end{document}